\documentclass[a4paper,reqno,12pt]{amsart}
\usepackage{enumerate}
\usepackage{amssymb}
\usepackage{ifthen}
\usepackage{graphicx}
\usepackage{color,xcolor}
\usepackage{mathrsfs}
\newtheorem{thm}{Theorem}
\newtheorem{theorem}{Theorem}

\newtheorem{lem}{Lemma}
\newtheorem{cor}{Corollary}
\newtheorem{defn}{Definition}
\newtheorem{example}{Example}

\newtheorem{claim}{Claim}

\newtheorem{conj}{Conjecture}

\newtheorem{prob}{Problem}

\newenvironment{rem}{%
\bigskip
\noindent \textsl{{\sl Remark. }}}{\bigskip}
\newenvironment{rems}{%
\bigskip
\noindent \textsl{{\sl Remarks. }}}{\bigskip}

\newenvironment{pf}[1][]{%
 \vskip 1mm
 \noindent
 \ifthenelse{\equal{#1}{}}%
  {{\slshape Proof. }}%
  {{\slshape #1.} }%
 }%
{\qed\bigskip}

\newcounter{alphabet}

\newenvironment{Lem}[1][]{\refstepcounter{alphabet}%
\bigskip%
\noindent%
{\bf Lemma \Alph{alphabet}}%
{\bf .} \itshape}{\vskip 8pt}

\DeclareMathOperator{\RE}{Re}

\newcommand{\ID}{{\mathbb D}}

\newcommand{\D}{{\mathbb D}}

\newcommand{\dist}{{\operatorname{dist}}}

\def\be{\begin{equation}}
\def\ee{\end{equation}}
\def\bes{\begin{equation*}}
\def\ees{\end{equation*}}

\newcommand{\bee}{\begin{enumerate}}
\newcommand{\eee}{\end{enumerate}}

\newcommand{\blem}{\begin{lem}}
\newcommand{\elem}{\end{lem}}
\newcommand{\bthm}{\begin{thm}}
\newcommand{\ethm}{\end{thm}}
\newcommand{\bcor}{\begin{cor}}
\newcommand{\ecor}{\end{cor}}
\newcommand{\beg}{\begin{example}}
\newcommand{\eeg}{\end{example}}
\newcommand{\begs}{\begin{examples}}
\newcommand{\eegs}{\end{examples}}
\newcommand{\bdefe}{\begin{defn}}
\newcommand{\edefe}{\end{defn}}
\newcommand{\bprob}{\begin{prob}}
\newcommand{\eprob}{\end{prob}}
\newcommand{\bques}{\begin{ques}}
\newcommand{\eques}{\end{ques}}
\newcommand{\bei}{\begin{itemize}}
\newcommand{\eei}{\end{itemize}}

\newcommand{\bde}{\begin{deter}}
\newcommand{\ede}{\end{deter}}
\newcommand{\bca}{\begin{case}}
\newcommand{\eca}{\end{case}}
\newcommand{\bcl}{\begin{claim}}
\newcommand{\ecl}{\end{claim}}
\newcommand{\bcon}{\begin{conj}}
\newcommand{\econ}{\end{conj}}
\newcommand{\bcons}{\begin{conjs}}
\newcommand{\econs}{\end{conjs}}
\newcommand{\bprop}{\begin{propo}}
\newcommand{\eprop}{\end{propo}}
\newcommand{\br}{\begin{rem}}
\newcommand{\er}{\end{rem}}
\newcommand{\brs}{\begin{rems}}
\newcommand{\ers}{\end{rems}}
\newcommand{\bo}{\begin{obser}}
\newcommand{\eo}{\end{obser}}
\newcommand{\bos}{\begin{obsers}}
\newcommand{\eos}{\end{obsers}}
\newcommand{\bpf}{\begin{pf}}
\newcommand{\epf}{\end{pf}}
\newcommand{\ba}{\begin{array}}
\newcommand{\ea}{\end{array}}
\newcommand{\beq}{\begin{eqnarray}}
\newcommand{\beqq}{\begin{eqnarray*}}
\newcommand{\eeq}{\end{eqnarray}}
\newcommand{\eeqq}{\end{eqnarray*}}

\newcommand{\ds}{\displaystyle}

\newcounter{minutes}
\divide\time by 60
\newcounter{hours}
\multiply\time by 60 \addtocounter{minutes}{-\time}

\begin{document}
\bibliographystyle{amsplain}
\title[Stable classes of harmonic mappings]
{Stable classes of harmonic mappings}

\def\thefootnote{}
\footnotetext{ \texttt{\tiny File:~\jobname .tex,
          printed: \number\day-\number\month-\number\year,
          \thehours.\ifnum\theminutes<10{0}\fi\theminutes}
} \makeatletter\def\thefootnote{\@arabic\c@footnote}\makeatother

\author[G. Liu]{Gang Liu}
\address{G. Liu, College of Mathematics and Statistics
 (Hunan Provincial Key Laboratory of Intelligent Information Processing and Application),
Hengyang Normal University, Hengyang,  Hunan 421002, China.}
\email{liugangmath@sina.cn}

\author[S. Ponnusamy]{Saminathan Ponnusamy
}
\address{S. Ponnusamy, Department of Mathematics,
Indian Institute of Technology Madras, Chennai-600 036, India;}
\address{Lomonosov Moscow State University, Moscow Center of Fundamental and Applied Mathematics, Moscow, Russia.
}
\email{samy@iitm.ac.in}
\author[V. V. Starkov]{Victor V. Starkov}
\address{V. V. Starkov, Department of Mathematics, University of Petrozavodsk,
ul. Lenina 33, 185910 Petrozavodsk, Russia.}
\email{vstarv@list.ru}

\subjclass[2010]{Primary: 30A10, 30B10, 30H05, 31A05, 30C62, 30C80; Secondary 30C35, 30C45}

\keywords{Harmonic mapping, univalent, convex, subordination, linear and affine invariant family,  Bohr inequality.
}

\begin{abstract}
Let $\mathcal{H}_0$ denote the set of all sense-preserving harmonic mappings $f=h+\overline{g}$ in the unit disk $\ID$, normalized with $h(0)=g(0)=g'(0)=0$ and $h'(0)=1$.
In this paper, we investigate some properties of certain subclasses  of $\mathcal{H}_0$,
including inclusion relations and stability analysis by precise examples, coefficient bounds, growth, covering and distortion theorems.
As applications, we build some Bohr inequalities for these subclasses by means of subordination.
Among these subclasses, six classes consist of functions $f=h+\overline{g}\in\mathcal{H}_0$ such that $h+\epsilon g$ is univalent (or convex) in $\D$
for each $|\epsilon|=1$ (or for some $|\epsilon|=1$, or for some $|\epsilon|\leq1$).
Simple analysis shows that if the function $f=h+\overline{g}$  belongs to a given class from  these six classes,
then the functions $h+\overline{\epsilon g}$ belong to corresponding class for all $|\epsilon|=1$.
We call these classes as stable classes.
\end{abstract}

\maketitle \pagestyle{myheadings}
\markboth{G. Liu, S. Ponnusamy and V. V. Starkov}{Stable classes of harmonic mappings}

\section{Introduction}  \label{sec1}

Let $\ID$ be the open unit disk, and $\mathbb{T}$ be the unit circle.
A complex-valued function $f=u+iv$ defined in $\ID$ is harmonic if  both $u$ and $v$ are real-valued harmonic functions in $\ID$.
Such a function can be represented as
$f=h+\overline{g}$, where $h$ and $g$ are analytic in $\ID$.
It is usually called {\it harmonic mapping}, and the function $h$ (resp. $g$) is called the {\it analytic part} (resp. {\it co-analytic part})  of $f$.
Lewy's theorem \cite{lew} states that $f=h+\overline{g}$ is locally univalent in $\mathbb{D}$ and  sense-preserving if and only if $h'(z)\neq0$
and its dilatation $\omega_f=g'/h'$  satisfies the condition $|\omega_f|<1$ in $\mathbb{D}$.

Let $\mathcal{S}_H$ be the family of all univalent, sense-preserving harmonic mappings $f=h+\overline{g}$ in $\ID$ normalized by $h(0)=g(0)=h'(0)-1=0$.
The class $\mathcal{S}_H$ is not compact whereas $\mathcal{S}_H^0:=\{f=h+\overline{g}\in\mathcal{S}_H:\,g'(0)=0\}$ is compact
(see \cite[p. 78]{dur2004}).
Note that the classical univalent class $\mathcal{S}$ is equivalent to the set $\{f=h+\overline{g}\in\mathcal{S}_H:~g(z)\equiv0\}$.
For every function $f\in\mathcal{S}$,  de Branges' theorem  \cite{deB} states that $|f^{(n)}(0)/n!|\leq n$ for all $n\geq2$,
and the equality is attained for the Koebe function $k(z)=z/(1-z)^2$ and its rotations.
This settles the Bieberbach conjecture.
For the class $\mathcal{S}_H^0$, Clunie and Sheil-Small proposed the following harmonic analogue of the Bieberbach conjecture.

\bcon{ \rm(\cite{CS})}\label{conjA}
For $f=h+\bar{g} \in \mathcal{S}_{H}^{0}$ with the series representation
\begin{equation}\label{equ1}
h(z)=z+\sum_{n=2}^\infty a_nz^n ~\mbox{ and }~g(z)=\sum_{ n=2}^\infty b_nz^n, ~z\in \ID,
\end{equation}
we have
$$|a_n|\leq A_n, ~|b_n|\leq B_n ~\mbox{ and }~ \big | |a_n|-|b_n| \big | \leq A_n-B_n=n
$$
for all $n\geq 2$, where
\beq\label{equ2}
A_{n}=\frac{(n+1)(2n+1)}{6}~\mbox{ and }~ B_{n}=\frac{(n-1)(2 n-1)}{6}.
\eeq
The bounds are attained for the harmonic Koebe function $K(z)$, defined by
\beq\label{equ3}
K(z)=\frac{z-\frac{1}{2} z^{2}+\frac{1}{6} z^{3}}{(1-z)^{3}}+\overline{\left(\frac{\frac{1}{2} z^{2}+\frac{1}{6} z^{3}}{(1-z)^{3}}\right)}.
\eeq
\econ

This conjecture remains an open problem, although it has been settled for a series of subclasses of $\mathcal{S}_H^0$,
for example,  the class of all starlike functions, close-to-convex, convex, typically real, or convex in one direction (see \cite{CS,she,WLZ}), respectively.
Ponnusamy and Sairam Kaliraj \cite{PS} certified that the result of coefficient bounds on $a_n$ and $b_n$ in Conjecture~A  
holds for the  class
$$\{f=h+\overline{g}\in\mathcal{S}_H^0:~h+e^{i\theta}g\in\mathcal{S}~\text{for~some}~\theta\in\mathbb{R}\}.
$$
A natural question is to ask whether or not there is a link on functions in the two families $\mathcal{S}_H^0$ and $\mathcal{S}$?
In view of this, they posed a conjecture that for every function $f=h+\overline{g}\in\mathcal{S}_H^0$,  there exists at least one  $\theta\in\mathbb{R}$ such that $h+e^{i\theta}g\in\mathcal{S}$.
Thus, such conjecture implies the estimates of $|a_n|$ and $|b_n|$ in Conjecture~A.  
Against the expectation, it was disproved in \cite{ABLS}.
Therefore, a  weaker conjecture can be considered.
That is, for every function $f=h+\overline{g}\in\mathcal{S}_H^0$,  there exists at least one
$\epsilon\in\overline{\mathbb{D}}$ such that $h+\epsilon g\in\mathcal{S}$.
It implies the sharp estimate on $|a_2|$ in Conjecture~A 
(see \cite{LP2019}).
For the families $\mathcal{C}_H^0$ and $\mathcal{C}$ of convex harmonic mappings in $\mathcal{S}_H^0$ and $\mathcal{S}$, similar conjecture was also listed in \cite{PS}.
Unfortunately, it is  also not true (see \cite{KPX}).
However, we do not know whether there exists at least one
$\epsilon\in\overline{\mathbb{D}}$ such that $h+\epsilon g\in\mathcal{C}$ when $f=h+\overline{g}\in\mathcal{C}_H^0$.

We are  not attempting to solve the above conjectures. However, along the line of the above weaker  conjectures, we will investigate some properties of the class of harmonic mappings $\mathcal{H}_0^*(\mathcal{S},\overline{\mathbb{D}})$ and related classes, 
where
$$\mathcal{H}_0^*(\mathcal{S},\overline{\mathbb{D}})=\left \{f=h+\overline{g}\in\mathcal{H}_0:\, h+\epsilon g\in\mathcal{S}~\text{
for~ some}  ~ \epsilon\in\overline{\mathbb{D}}\right \}.
$$
Here $\mathcal{H}_0$ is the family of sense-preserving harmonic mappings $f=h+\overline{g}$ in $\ID$, normalized with $h(0)=g(0)=g'(0)=0$ and $h'(0)=1$.

The paper is organized as follows.
Before investigating some properties of $\mathcal{H}_0^*(\mathcal{S},\overline{\mathbb{D}})$, we will introduce its five proper subclasses in Section \ref{sec2}.
The class $\mathcal{H}_0^*(\mathcal{S},\overline{\mathbb{D}})$ and its five proper subclasses 
will be called stable classes of harmonic mappings,
and they will be divided into two groups with the same number.
Since the two groups have similar structures in common, 
 below we only 
 introduce our results concerning one group related to $\mathcal{S}$.
The three classes in the first group have corresponding functions with different stabilities from the analysis of their strict inclusion relations.
In Section \ref{sec3}, we will show that the former two coefficient bounds in Conjecture~A 
hold for functions in $\mathcal{H}_0^*(\mathcal{S},\overline{\mathbb{D}})$.
This means that if $\mathcal{S}_H^0\subseteq\mathcal{H}_0^*(\mathcal{S},\overline{\mathbb{D}})$ was true, then we can solve Conjecture~A 
partly.
However, it seems very difficult to prove or disprove it now.
On the basis of the fact that the class $\mathcal{A}(\mathcal{F})$ is an affine and linear invariant family if $\mathcal{F}$ is chosen from anyone from the six stable subclasses of $\mathcal{H}_0^*(\mathcal{S},\overline{\mathbb{D}})$, where
$$\mathcal{A}(\mathcal{F}) :=\left \{F=f+\overline{b_1f}:\,b_1\in\ID ~\mbox{ and }~ f\in\mathcal{F}\subseteq\mathcal{H}_0 \right \}.
$$
In  Section \ref{sec4}, we can obtain some growth, covering and distortion theorems for functions in certain subclasses of $\mathcal{A}(\mathcal{H}_0^*(\mathcal{S},\overline{\mathbb{D}}))$.
As applications, in Section \ref{sec5}, we will build several Bohr inequalities for the class of functions subordinate to functions in $\mathcal{H}_0^*(\mathcal{S},\overline{\mathbb{D}})$ and its subclasses.
It is worth mentioning that for functions in two different stable subclasses of $\mathcal{H}_0^*(\mathcal{S},\overline{\mathbb{D}})$, although their stabilities are different,
they  share some common properties and are investigated in the last three sections.

\section{Stability analysis of $\mathcal{H}_0^*(\mathcal{S},\overline{\mathbb{D}})$ and its subclasses} \label{sec2}

The aim of this section is three-fold.
One is to introduce  a notion of stable class of harmonic mappings since the class $\mathcal{H}_0^*(\mathcal{S},\overline{\mathbb{D}})$ and its five proper subclasses  have a common algebraic characterization.
The second is to exhibit different stabilities for two groups that each has three stable subclasses of $\mathcal{H}_0^*(\mathcal{S},\overline{\mathbb{D}})$,
which are obtained from their strict inclusion relations through precise examples.
The last is to list certain inclusion relations between  the two non-stable classes $\mathcal{S}_H^0$ and $\mathcal{C}_H^0$, and six stable subclasses of $\mathcal{H}_0^*(\mathcal{S},\overline{\mathbb{D}})$.
Below we will frequently use the following two notations,  for a given $f=h+\overline{g}\in\mathcal{H}:=\mathcal{A}(\mathcal{H}_0)$,
$$f_\epsilon=h+\overline{\epsilon g}\quad \text {and} \quad f^\epsilon=h+\epsilon g ~(\epsilon\in\overline{\mathbb{D}}).
$$
Here, the five subclasses of $\mathcal{H}_0^*(\mathcal{S},\overline{\mathbb{D}})$ are defined as follows:
$$\mathcal{H}_0^*(\mathcal{C},\overline{\mathbb{D}})=\left \{f\in\mathcal{H}_0:\, \mbox{$f^\epsilon\in\mathcal{C}$
~ for some  $\epsilon\in\overline{\mathbb{D}}$}\right \},$$
$$\mathcal{H}_0^*(\mathcal{S},\mathbb{T})~(\mbox{resp. $\mathcal{H}_0^*(\mathcal{C},\mathbb{T})$})=\left \{f\in\mathcal{H}_0:\, \mbox{$f^\epsilon\in\mathcal{S}$
~(\mbox{resp. $\mathcal{C}$}) for some  $\epsilon\in\mathbb{T}$}\right \}$$
and
$$\mathcal{H}_0(\mathcal{S},\mathbb{T})~(\mbox{resp. $\mathcal{H}_0(\mathcal{C},\mathbb{T})$})=\left \{f\in\mathcal{H}_0:\, \mbox{$f^\epsilon\in\mathcal{S}$
~(\mbox{resp. $\mathcal{C}$}) for every  $\epsilon\in\mathbb{T}$}\right \}.$$
Clearly,  $\mathcal{S}\subset\mathcal{H}_0(\mathcal{S},\mathbb{T})$ and $\mathcal{C}\subset\mathcal{H}_0(\mathcal{C},\mathbb{T})$ by considering the function $z+\overline{z^2/2}$.
It follows from $\mathcal{C}\subset\mathcal{S}$  that  $\mathcal{H}_0(\mathcal{C},\mathbb{T})\subset\mathcal{H}_0(\mathcal{S},\mathbb{T})$, $\mathcal{H}_0^*(\mathcal{C},\mathbb{T})\subset\mathcal{H}_0^*(\mathcal{S},\mathbb{T})$ and  $\mathcal{H}_0^*(\mathcal{C},\overline{\mathbb{D}})\subset\mathcal{H}_0^*(\mathcal{S},\overline{\mathbb{D}})$.
Note that $\mathcal{H}_0(\mathcal{S},\mathbb{T})\subseteq\mathcal{H}_0^*(\mathcal{S},\mathbb{T})\subseteq\mathcal{H}_0^*(\mathcal{S},\overline{\mathbb{D}})$ and $\mathcal{H}_0(\mathcal{C},\mathbb{T})\subseteq\mathcal{H}_0^*(\mathcal{C},\mathbb{T})\subseteq\mathcal{H}_0^*(\mathcal{C},\overline{\mathbb{D}})$ from their definitions.

\subsection{Stable classes of harmonic mappings}

The family $\mathcal{F}\subseteq\mathcal{H}$ is called a stable class of harmonic mappings (in what follows, for simplicity,
we may call that $\mathcal{F}$ is stable), if $f\in\mathcal{F}$, then $f_\epsilon\in\mathcal{F}$ for every $\epsilon\in\mathbb{T}$.
In other words, $f_{\epsilon_1}\in\mathcal{F}$ is equivalent to $f_{\epsilon_2}\in\mathcal{F}$ for any pairs $\epsilon_1,~\epsilon_2\in\mathbb{T}$ with $\epsilon_1\neq\epsilon_2$.

The terminology ``stable" is adapted from the notion of stable mapping investigated by Chuaqui and  Hern\'{a}ndez \cite{CH}, and  Hern\'{a}ndez and Mart\'{\i}n \cite{HM}.
A sense-preserving harmonic mapping $f=h+\overline{g}$  in $\mathbb{D}$ is called stable harmonic univalent or SHU
(resp.  stable harmonic convex or SHC)  if $f_\epsilon$ is univalent (resp. convex)  in $\mathbb{D}$  for every $\epsilon\in\mathbb{T}$.
It has been certified that for all $\epsilon\in\mathbb{T}$, the harmonic mappings $f_\epsilon$ are univalent (resp. convex)
if and only if the analytic functions $f^\epsilon$ are univalent (resp. convex) for all such $\epsilon$ (see \cite{HM}).
This means that if $f$ is SHU (resp. SHC), then each function $f_\epsilon$ is SHU (resp. SHC) for every $\epsilon\in\mathbb{T}$.
Therefore, the set $\mathcal{H}_0(\mathcal{S},\mathbb{T})$ (resp. $\mathcal{H}_0(\mathcal{C},\mathbb{T})$) consists of all functions that are SHU (resp. SHC) in the class $\mathcal{H}_0^*(\mathcal{S},\overline{\mathbb{D}})$.

Let us show some simple properties of stable classes of harmonic mappings.
Assume that the classes $\mathcal{F}_1$ and $\mathcal{F}_2$ are stable.
It is clearly that both $\mathcal{F}_1\cup\mathcal{F}_2$, $\mathcal{F}_1\cap\mathcal{F}_2\, (\neq\emptyset)$ and $\mathcal{F}_1\setminus\mathcal{F}_2\, (\neq\emptyset)$ are stable.
Furthermore, if we find a function $f\in\mathcal{F}_1\setminus\mathcal{F}_2$, then we know that $f_\epsilon\in\mathcal{F}_1\setminus\mathcal{F}_2$ for every $\epsilon\in\mathbb{T}$.
This 
means that if $f=h+\overline{g}\in\mathcal{F}_1\setminus\mathcal{F}_2$ and $g(z)\not\equiv0$ in $\ID$, then   there are infinitely many functions in $\mathcal{F}_1\setminus\mathcal{F}_2$.

If the family $\mathcal{F}_0\subseteq\mathcal{H}_0$ is stable, then the family $\mathcal{A}(\mathcal{F}_0)$ is also stable.
For this, if $F\in\mathcal{A}(\mathcal{F}_0)$, then there  exists a constant $b_1\in\D$ and a function $f=h+\overline{g}\in\mathcal{F}_0$  such that $F=f+\overline{b_1f}=h+\overline{b_1}g+\overline{b_1h+g}$.
Since $f\in\mathcal{F}_0$ implies $f_\epsilon\in\mathcal{F}_0$ for every $\epsilon\in\mathbb{T}$, a simple computation shows that
$$F_\epsilon=h+\overline{b_1}g+\overline{\epsilon(b_1h+g)}=h+\overline{\epsilon g}+\overline{(b_1\epsilon)(h+\overline{\epsilon g})}=f_\epsilon+\overline{(b_1\epsilon)f_\epsilon}\in\mathcal{A}(\mathcal{F}_0)
$$
for every $\epsilon\in\mathbb{T}$.
It is easy to see that the class $\mathcal{F}_0$ is stable if  $\mathcal{F}_0$ is chosen by anyone of the seven families:
$\mathcal{H}_0$, $\mathcal{H}_0^*(\mathcal{S},\overline{\mathbb{D}})$, $\mathcal{H}_0^*(\mathcal{C},\overline{\mathbb{D}})$, $\mathcal{H}_0^*(\mathcal{S},\mathbb{T})$, $\mathcal{H}_0^*(\mathcal{C},\mathbb{T})$, $\mathcal{H}_0(\mathcal{S},\mathbb{T})$ and $\mathcal{H}_0(\mathcal{C},\mathbb{T})$. So does  with its corresponding $\mathcal{A}(\mathcal{F}_0)$.

\subsection{Stability analysis of $\mathcal{H}_0^*(\mathcal{S},\overline{\mathbb{D}})$ and its subclasses}
The stability of a stable class of harmonic mappings is similar to search for the relationship
between $f_\epsilon$ and $f^\epsilon$ when the function $f$ is  restricted to be in SHU or SHC, as investigated in \cite{HM}.
Let us first describe the stability of the two classes $\mathcal{H}_0(\mathcal{S},\mathbb{T})$ and $\mathcal{H}_0(\mathcal{C},\mathbb{T})$.
It follows from \cite[Theorem~7.1]{HM} that if $f\in\mathcal{H}_0(\mathcal{S},\mathbb{T})$,
i.e., $f$ is a stable univalent  mapping from $\mathcal{H}_0^*(\mathcal{S},\overline{\mathbb{D}})$,
then $f_\epsilon\in\mathcal{H}_0(\mathcal{S},\mathbb{T})$ and $f^\epsilon\in\mathcal{S}$ for every $\epsilon\in\overline{\D}$.
Similarly, it follows from \cite[Corollary 3.2]{HM} that if $f\in\mathcal{H}_0(\mathcal{C},\mathbb{T})$,
i.e., $f$ is a stable convex mapping from $\mathcal{H}_0^*(\mathcal{S},\overline{\mathbb{D}})$,
then $f_\epsilon\in\mathcal{H}_0(\mathcal{C},\mathbb{T})$ and $f^\epsilon\in\mathcal{C}$ for every $\epsilon\in\overline{\D}$.
Thus, the class $\mathcal{H}_0(\mathcal{S},\mathbb{T})$ (resp. $\mathcal{H}_0(\mathcal{C},\mathbb{T})$) can be rewritten as
\beq\label{equ4}
\mathcal{H}_0(\mathcal{S},\mathbb{T})~(\mbox{resp. $\mathcal{H}_0(\mathcal{C},\mathbb{T})$})=\left \{f\in\mathcal{H}_0:\, \mbox{$f^\epsilon\in\mathcal{S}$
~(\mbox{resp. $\mathcal{C}$}) for every  $\epsilon\in\overline{\mathbb{D}}$}\right \}.
\eeq

 Comparing with  the stability of the two classes $\mathcal{H}_0(\mathcal{S},\mathbb{T})$ and $\mathcal{H}_0(\mathcal{C},\mathbb{T})$,
the remaining four stable subclasses of $\mathcal{H}_0^*(\mathcal{S},\overline{\mathbb{D}})$ mentioned  earlier have weaker stability.
 Indeed, based on the two pairs of strict inclusion relations by  specific examples, we will exhibit their differences on stability by two tables.

We claim that $\mathcal{H}_0(\mathcal{S},\mathbb{T})\subset\mathcal{H}_0^*(\mathcal{S},\mathbb{T})\subset\mathcal{H}_0^*(\mathcal{S},\overline{\mathbb{D}})$.
Recall that the harmonic Koebe function $K=h+\overline{g}$ defined in \eqref{equ3}, where
\beqq
h(z)=\frac{z-\frac{1}{2} z^{2}+\frac{1}{6} z^{3}}{(1-z)^{3}}=z+\sum_{n=2}^{\infty} A_{n} z^{n}\quad \text { and }
\quad g(z)=\frac{\frac{1}{2} z^{2}+\frac{1}{6} z^{3}}{(1-z)^{3}}=\sum_{n=2}^{\infty} B_{n} z^{n}.
\eeqq
Here $A_n$ and $B_n$ are defined by \eqref{equ2}.
Direct computation gives the following power series representation for $K^\epsilon =h+\epsilon g$:
$$
K^\epsilon(z)=z+\sum_{n=2}^{\infty} \varphi_n(\epsilon) z^{n},
$$
where $\epsilon\in\overline{\D}$ and
$$
\varphi_n(\epsilon)=\frac{1}{6}\big [(2 n^{2}+1)(1+\epsilon)+3 n(1-\epsilon)] ~\mbox{ for all $ n \geq 2$}.
$$
Then $K^\epsilon\in\mathcal{S}$ only for $\epsilon=-1$ in which case $K^\epsilon$ reduces to the analytic Koebe function  $k(z)=z/(1-z)^2$,
or else it contradicts  de Branges' coefficient inequality when $n$ is sufficient large.
This implies that $K\in\mathcal{H}_0^*(\mathcal{S},\mathbb{T})$ but $K\not\in\mathcal{H}_0(\mathcal{S},\mathbb{T})$.

Recall that $\mathcal{H}_0(\mathcal{S},\mathbb{T})\subseteq\mathcal{H}_0^*(\mathcal{S},\mathbb{T})\subseteq\mathcal{H}_0^*(\mathcal{S},\overline{\mathbb{D}})$.
Next it suffices to find a function in $\mathcal{H}_0^*(\mathcal{S},\overline{\mathbb{D}})\backslash\mathcal{H}_0^*(\mathcal{S},\mathbb{T})$.
Let $M=h+\overline{g}$, where $h(z)=k(z)=z/(1-z)^2$ and
$$g(z)=\frac{-z(24-35z+8z^2+z^3)}{6(1-z)^2}-4\log(1-z)=-\frac{1}{6}z^2+\cdots+\frac{5}{21}z^7+\cdots.
$$
Clearly, $g(0)=g'(0)=0$. Moreover, a simple calculation shows that
$$g'(z)=\frac{z(-1+z+z^2)(1+z)}{3(1-z)^3}
$$
and thus, $M\in\mathcal{H}_0$.
We will prove that $M^\epsilon\in\mathcal{S}~(\epsilon\in\overline{\mathbb{D}})$ only if $\epsilon=0$.
To do this we first set
$$M^\epsilon(z)= h(z)+ \epsilon g(z) = z+\cdots+a_7(\epsilon)z^7+\cdots
$$
and consider the following two sets
$$  R(\mathbb{D})=\overline{\mathbb{D}}\cap\{z=x+iy:\, x\geq0\}
~\mbox{ and }~ L(\mathbb{D})=\overline{\mathbb{D}}\cap\{z=x+iy:\, x\leq0\}.
$$
If $\epsilon\in R(\mathbb{D})\setminus \{0\}$, then we have $|a_7(\epsilon)|=|7+(5/21)\epsilon|>7$, which implies $M^\epsilon\not\in\mathcal{S}$ when $\epsilon\in R(\mathbb{D})\setminus \{0\}$ by de Branges' theorem.
If $\epsilon\in L(\mathbb{D})\setminus \{0\}$, then for $0<r<(\sqrt{5}-1)/2$, we get
$$|(M^\epsilon)'(r)|=\left|1-\frac{1}{3}r(1-r-r^2)\epsilon\right|\frac{1+r}{(1-r)^3}>\frac{1+r}{(1-r)^3}
$$
since  $1-r-r^2>0$ whenever $0<r<(\sqrt{5}-1)/2$.
This shows that $M^\epsilon\not\in\mathcal{S}$ whenever $\epsilon\in L(\mathbb{D})\setminus \{0\}$ by the classical distortion theorem for functions in the class $\mathcal{S}$ (see \cite[Theorem~1.6]{pom}).
Therefore, we get $M^\epsilon =h+\epsilon g\not\in\mathcal{S}$ for all $\epsilon\in\overline{\mathbb{D}}\setminus\{0\}$, because of $R(\mathbb{D})\cup L(\mathbb{D})=\overline{\D}$. Note that $M^0(z)=k(z)\in\mathcal{S}$.
It follows that $M\in\mathcal{H}_0^*(\mathcal{S},\overline{\mathbb{D}})$ but $M\not\in\mathcal{H}_0^*(\mathcal{S},\mathbb{T})$.
Thus, the desired strict inclusions $\mathcal{H}_0(\mathcal{S},\mathbb{T})\subset\mathcal{H}_0^*(\mathcal{S},\mathbb{T})\subset\mathcal{H}_0^*(\mathcal{S},\overline{\mathbb{D}})$ follow.
Further analysis gives the following conclusions in Table \ref{tab1} which exhibit their different stabilities.
Note that the parameter $\epsilon\in\overline{\mathbb{D}}$.

\begin{table}[htb]
\begin{center}
\begin{tabular}{|c|c|c|c|}
\hline  1 & $\forall f\in\mathcal{H}_0(\mathcal{S},\mathbb{T})$  & $\exists K\in\mathcal{H}_0^*(\mathcal{S},\mathbb{T})\setminus\mathcal{H}_0(\mathcal{S},\mathbb{T})$ & $\exists M\in\mathcal{H}_0^*(\mathcal{S},\overline{\mathbb{D}})\setminus\mathcal{H}_0^*(\mathcal{S},\mathbb{T})$ \\
\hline  2  &$f_\epsilon\in\mathcal{H}_0(\mathcal{S},\mathbb{T})$ iff $\epsilon\in\overline{\mathbb{D}}$ & $K_\epsilon\in\mathcal{H}_0^*(\mathcal{S},\mathbb{T})$ iff $\epsilon\in\mathbb{T}$ & $M_\epsilon\in\mathcal{H}_0^*(\mathcal{S},\overline{\mathbb{D}})$ iff $\epsilon\in\overline{\mathbb{D}}$ \\
\hline  3  & $f^\epsilon\in\mathcal{S}$ iff $\epsilon\in\overline{\mathbb{D}}$  & $K^\epsilon\in\mathcal{S}$ iff $\epsilon=-1 ~(\in\mathbb{T})$ & $M^\epsilon\in\mathcal{S}$ iff $\epsilon=0~ (\in\mathbb{D})$  \\
\hline
\end{tabular}
\end{center}
\caption{Stability analysis of $\mathcal{H}_0(\mathcal{S},\mathbb{T})$, $\mathcal{H}_0^*(\mathcal{S},\mathbb{T})$ and $\mathcal{H}_0^*(\mathcal{S},\overline{\mathbb{D}})$}
\label{tab1}
\end{table}

We also state that $\mathcal{H}_0(\mathcal{C},\mathbb{T})\subset\mathcal{H}_0^*(\mathcal{C},\mathbb{T})\subset\mathcal{H}_0^*(\mathcal{C},\overline{\mathbb{D}})$.
Consider the harmonic half-plane mapping $L=h+\overline{g}$ investigated in \cite{CS}, where
\beqq
h(z)=\frac{z- z^{2}/2}{(1-z)^2}=z+\sum_{n=2}^{\infty} \frac{n+1}{2} z^{n}\quad \text { and }
\quad g(z)=\frac{-z^{2}/2}{(1-z)^2}=-\sum_{n=2}^{\infty} \frac{n-1}{2} z^{n}.
\eeqq
Now, a straightforward computation yields
$$L^\epsilon(z)=z+\sum_{n=2}^{\infty} \phi_n(\epsilon) z^{n},
$$
where $L^\epsilon =h+\epsilon g$ and
$$ \phi_n(\epsilon)=\frac{1}{2}\left[n(1-\epsilon)+(1+\epsilon)\right] ~\mbox{ for all $ n \geq 2$}.
$$
We observe $|\phi_n(\epsilon)|>1$ when $\epsilon\in\overline{\mathbb{D}}\setminus\{1\}$ and $n$ is sufficient large.
Thus, we know that $L^\epsilon =h+\epsilon g\in\mathcal{C}$  only for $\epsilon=1$ in which case it reduces to the
analytic half-plane function $l(z)=z/(1-z)$. This shows that $L\in\mathcal{H}_0^*(\mathcal{C},\mathbb{T})$ but $L\not\in\mathcal{H}_0(\mathcal{C},\mathbb{T})$.

The desired claim easily follows from  the observation that $\mathcal{H}_0^*(\mathcal{C},\overline{\mathbb{D}})\backslash\mathcal{H}_0^*(\mathcal{C},\mathbb{T})\neq\emptyset$, since $\mathcal{H}_0(\mathcal{C},\mathbb{T})\subseteq\mathcal{H}_0^*(\mathcal{C},\mathbb{T})\subseteq\mathcal{H}_0^*(\mathcal{C},\overline{\mathbb{D}})$.
Let $V_\alpha=h+\overline{g_\alpha}$, where $h(z)=l(z)=z/(1-z)$ and
$$g_{\alpha}(z)=\frac{\alpha z^{n}}{1-z},\quad 0<|\alpha|\leq\frac{1}{2n-1}\quad \text{and}\quad n\in\mathbb{N}\backslash \{1,2\}.
$$
The family of functions $V_\alpha\in\mathcal{H}_0$ was introduced in \cite[Section~3.1]{PS}.
Below we will use the result that if $f\in\mathcal{C}$, then $\RE\frac{f(z)}{z}>\frac{1}{2}$ for all $z\in\D$ (see \cite[Corollary~2]{suf}).
Fix $n$ and $\alpha$. For any given $\epsilon\in\overline{\mathbb{D}}\setminus\{0\}$,
let $z_0=re^{i\varphi}~(0<r<1)$ and $\varphi=(\pi-\arg\alpha-\arg\epsilon)/(n-1)$.
We can deduce that
$$\RE\frac{(V_\alpha)^\epsilon(z_0)}{z_0}=\frac{(1-|\alpha\epsilon|r^{n-1})(1-r\cos\varphi)}{1+r^2-2r\cos\varphi}\rightarrow\frac{1-|\alpha\epsilon|}{2}\quad(r\rightarrow1^-).
$$
Since $\frac{1-|\alpha\epsilon|}{2}<\frac{1}{2}$, we get $(V_\alpha)^\epsilon\not\in\mathcal{C}$ by continuity.
Note that $(V_\alpha)^0(z)=l(z)\in\mathcal{C}$. Hence, $V_\alpha\in\mathcal{H}_0^*(\mathcal{C},\overline{\mathbb{D}})\setminus\mathcal{H}_0^*(\mathcal{C},\mathbb{T})$.
Therefore, we can get the strict inclusion relations $\mathcal{H}_0(\mathcal{C},\mathbb{T})\subset\mathcal{H}_0^*(\mathcal{C},\mathbb{T})\subset\mathcal{H}_0^*(\mathcal{C},\overline{\mathbb{D}})$.
Furthermore, we can obtain the conclusions in  Table \ref{tab2}  so that each class in the three classes has distinct stability.
Note that the parameter $\epsilon\in\overline{\mathbb{D}}$.

\begin{table}[htb]
\begin{center}
 \begin{tabular}{|c|c|c|c|}
\hline   1 &$\forall f\in\mathcal{H}_0(\mathcal{C},\mathbb{T})$  & $\exists L\in\mathcal{H}_0^*(\mathcal{C},\mathbb{T})\setminus\mathcal{H}_0(\mathcal{C},\mathbb{T})$ & $\exists V_\alpha\in\mathcal{H}_0^*(\mathcal{C},\overline{\mathbb{D}})\setminus\mathcal{H}_0^*(\mathcal{C},\mathbb{T})$ \\
\hline  2 & $f_\epsilon\in\mathcal{H}_0(\mathcal{C},\mathbb{T})$ iff $\epsilon\in\overline{\mathbb{D}}$ & $L_\epsilon\in\mathcal{H}_0^*(\mathcal{C},\mathbb{T})$ iff $\epsilon\in\mathbb{T}$ & $(V_\alpha)_\epsilon\in\mathcal{H}_0^*(\mathcal{C},\overline{\mathbb{D}})$ iff $\epsilon\in\overline{\mathbb{D}}$ \\
\hline  3 & $f^\epsilon\in\mathcal{C}$ iff $\epsilon\in\overline{\mathbb{D}}$  & $L^\epsilon\in\mathcal{C}$ iff $\epsilon=-1 (\in\mathbb{T})$ & $(V_\alpha)^\epsilon\in\mathcal{C}$ iff $\epsilon=0 (\in\mathbb{D})$  \\
\hline
\end{tabular}
\end{center}
\caption{Stability analysis of $\mathcal{H}_0(\mathcal{C},\mathbb{T})$, $\mathcal{H}_0^*(\mathcal{C},\mathbb{T})$ and $\mathcal{H}_0^*(\mathcal{C},\overline{\mathbb{D}})$}
\label{tab2}
\end{table}

\subsection{Relationships with two non-stable classes} \label{rem1}
We claim that neither $\mathcal{S}_H^0$ nor $\mathcal{C}_H^0$ is stable.
Recall that $K\in\mathcal{H}_0^*(\mathcal{S},\mathbb{T})\setminus\mathcal{H}_0(\mathcal{S},\mathbb{T})$.
This implies that there exists at least one $\epsilon_0\in\mathbb{T}\setminus\{1\}$ such that $K_{\epsilon_0}$ is not univalent in $\D$.
It follows from $K\in\mathcal{S}_H^0$  (see \cite{CS}) that $\mathcal{S}_H^0$ is not stable.
For non-stability of the class $\mathcal{C}_H^0$, similar to the above analysis, it follows from the fact that $L\in\mathcal{C}_H^0$ (see \cite{CS})
and $L\in\mathcal{H}_0^*(\mathcal{C},\mathbb{T})\setminus\mathcal{H}_0(\mathcal{C},\mathbb{T})$.

Next we consider some relationships 
among $\mathcal{S}_H^0$, $\mathcal{C}_H^0$ and six stable subclasses of  $\mathcal{H}_0^*(\mathcal{S},\overline{\mathbb{D}})$.
It follows from \cite[Theorem~5.7]{CS} and \cite[Theorem~5.17]{CS} that $\mathcal{C}_H^0\subseteq\mathcal{H}_0(\mathcal{S},\mathbb{T})$ and $\mathcal{H}_0^*(\mathcal{C},\overline{\mathbb{D}})\subseteq\mathcal{S}_H^0$, respectively.
Note that $k\in\mathcal{H}_0(\mathcal{S},\mathbb{T})\setminus\mathcal{C}_H^0$ and  $K\in\mathcal{S}_H^0\setminus\mathcal{H}_0^*(\mathcal{C},\overline{\mathbb{D}})$.
Thus, we have
$$\mathcal{H}_0(\mathcal{C},\mathbb{T})\subset\mathcal{C}_H^0\subset\mathcal{H}_0(\mathcal{S},\mathbb{T})\subset\mathcal{S}_H^0\quad\text{and}\quad
\mathcal{H}_0^*(\mathcal{C},\overline{\mathbb{D}})\subset\mathcal{S}_H^0.
$$

However, we have $\mathcal{H}_0^*(\mathcal{S},\mathbb{T})\not\subset\mathcal{S}_H^0$, $\mathcal{H}_0^*(\mathcal{S},\overline{\mathbb{D}})\not\subset\mathcal{S}_H^0$ and $\mathcal{S}_H^0\not\subset\mathcal{H}_0^*(\mathcal{S},\mathbb{T})$.
From the above analysis, we know that $K_{\epsilon_0}\in\mathcal{H}_0^*(\mathcal{S},\mathbb{T})\setminus\mathcal{S}_H^0$  for some $\epsilon_0\in\mathbb{T}\setminus\{1\}$ , which in turn implies  that $\mathcal{H}_0^*(\mathcal{S},\mathbb{T})\not\subset\mathcal{S}_H^0$,
and thus $\mathcal{H}_0^*(\mathcal{S},\overline{\mathbb{D}})\not\subset\mathcal{S}_H^0$.
For $\mathcal{S}_H^0\not\subset\mathcal{H}_0^*(\mathcal{S},\mathbb{T})$, it is enough to find a function $f$ such that $f\in\mathcal{S}_H^0\setminus\mathcal{H}_0^*(\mathcal{S},\mathbb{T})$.
In other words, the function $f=h+\overline{g}\in\mathcal{S}_H^0$ is satisfied with the condition $h+\epsilon g\not\in\mathcal{S}$ for all $\epsilon\in\mathbb{T}$.
The construction of such  function  is not easy. For a precise example with technical details, we refer to the recent article \cite{ABLS}.


Similarly, we can get $\mathcal{H}_0^*(\mathcal{C},\mathbb{T})\not\subset\mathcal{C}_H^0$, $\mathcal{H}_0^*(\mathcal{C},\overline{\mathbb{D}})\not\subset\mathcal{C}_H^0$ and $\mathcal{C}_H^0\not\subset\mathcal{H}_0^*(\mathcal{C},\mathbb{T})$.
The desired conclusions follow from the fact that $L_{\epsilon_0}\in\mathcal{H}_0^*(\mathcal{C},\mathbb{T})\setminus\mathcal{C}_H^0$ for some $\epsilon_0\in\mathbb{T}\setminus\{1\}$ and the example constructed in \cite{KPX}.

It is natural to ask: {\em what are the precise inclusion relations among the following seven families?:} $\mathcal{S}_H^0$, $\mathcal{C}_H^0$,
$\mathcal{H}_0(\mathcal{S},\mathbb{T})$,  $\mathcal{H}_0^*(\mathcal{S},\mathbb{T})$, $\mathcal{H}_0^*(\mathcal{S},\overline{\mathbb{D}})$,  $\mathcal{H}_0^*(\mathcal{C},\mathbb{T})$ and $\mathcal{H}_0^*(\mathcal{C},\overline{\mathbb{D}})$.
In particular,  we ask whether the following inclusions hold: $\mathcal{S}_H^0\subset\mathcal{H}_0^*(\mathcal{S},\overline{\mathbb{D}})$ and $\mathcal{C}_H^0\subset\mathcal{H}_0^*(\mathcal{C},\overline{\mathbb{D}})$.
If $\mathcal{S}_H^0\subset\mathcal{H}_0^*(\mathcal{S},\overline{\mathbb{D}})$, 
then we can solve conjecture~A 
partly (see Theorem \ref{thm1}).
To the best of our knowledge, we cannot provide an answer now.

\section{Coefficient bounds for functions in the subclasses of $\mathcal{H}_0^*(\mathcal{S},\overline{\mathbb{D}})$}  \label{sec3}

In this section we will find bounds on coefficients  of functions in the six stable subclasses of $\mathcal{H}_0^*(\mathcal{S},\overline{\mathbb{D}})$.
The proof of the following result is similar to \cite[Theorem~1]{PS}.
For completeness, we present a proof with slight modifications.

\begin{theorem} \label{thm1} Suppose that the function $f=h+\bar{g}$  belongs to $\mathcal{H}_0^*(\mathcal{S},\overline{\mathbb{D}})$ or $\mathcal{H}_0^*(\mathcal{S},\mathbb{T})$, where $h$ and $g$ are given by \eqref{equ1}.
Then $|a_n|\leq A_n$ and $|b_n|\leq B_n$ for all $n\geq 2$, where $A_n$ and $B_n$ are defined by \eqref{equ2}.
These bounds are attained for the harmonic Koebe function $K$.
\end{theorem}

\bpf
Suppose that $f=h+\bar{g} \in \mathcal{H}_0^*(\mathcal{S},\overline{\mathbb{D}}).$
Then we have $f^\epsilon(z)=z+\sum_{n=2}^{\infty} \epsilon_{n} z^{n} \in \mathcal{S}$ for some $\epsilon\in\overline{\D}$,
where $f^\epsilon =h+ \epsilon g$.
By de Branges' theorem, $|\epsilon_{n}| \leq n$ for all $n \geq 2$.
Since $f$ is sense-preserving in $\mathbb{D},$ there exists an analytic function $\omega(z)$ in $\mathbb{D}$ such that $\omega(0)=0$ and $|\omega(z)|=\left|g'(z) / h'(z)\right|<1$ for all $z \in \mathbb{D}$ from which we easily obtain that
$$ \varphi'(z)=h'(z)(1+\epsilon \omega(z))  \quad \text { and } \quad g'(z)=\frac{\varphi'(z) \omega(z)}{1+\epsilon \omega(z)},
$$
where for convenience, we write $\varphi=f^\epsilon$.
We observe that the Taylor coefficients of the function $\varphi'(z)$ (resp. $\frac{\omega(z)}{1+\epsilon\omega(z)}$) are dominated in modulus by the corresponding coefficients of the function
$\frac{1+z}{(1-z)^3}$ (resp. $\frac{z}{1-z}$).
Therefore, the Taylor coefficients $nb_n$ of $g'(z)$ are dominated by those of
$$\frac{1+z}{(1-z)^3}\cdot\frac{z}{1-z}=\frac{1}{6}\sum_{n=2}^\infty(n-1)n(2n-1)z^{n-1},
$$
which implies $|b_n|\leq B_n =(n-1)(2 n-1)/6$ for $n\geq2$.
From the definition of $\varphi(z),$ we have $h(z)=\varphi(z)-\epsilon g(z) .$
Therefore, one has
$$|a_{n}|=\left|\epsilon_{n}-\epsilon b_{n}\right| \leq\left|\epsilon_{n}\right|+\left|b_{n}\right| \leq \frac{(n+1)(2 n+1)}{6} =A_n \quad \text { for } n \geq 2.
$$
Note that $\mathcal{H}_0^*(\mathcal{S},\mathbb{T})\subset\mathcal{H}_0^*(\mathcal{S},\overline{\mathbb{D}})$ and $K\in\mathcal{H}_0^*(\mathcal{S},\mathbb{T})$.
We complete the proof.
\epf

It is easy to get the following result.

\bcor \label{cor1}
 Suppose that the function $f(z)=h(z)+\overline{g(z)}=z+\sum_{n=2}^{\infty} a_n z^n+\overline{\sum_{n=1}^{\infty} b_n z^n}$
belongs to $\mathcal{A}(\mathcal{H}_0^*(\mathcal{S},\overline{\mathbb{D}}))$ or $\mathcal{A}(\mathcal{H}_0^*(\mathcal{S},\mathbb{T}))$.
Then
$$\left|a_{n}\right|<\frac{1}{3}\left(2 n^{2}+1\right) ~\mbox{ and }~ \left|b_{n}\right|<\frac{1}{3}\left(2 n^{2}+1\right)
~\mbox{ for all $n \geq 2$.}
$$
\ecor

Next we consider coefficient bounds for functions in classes $\mathcal{H}_0^*(\mathcal{C},\overline{\mathbb{D}})$ and $\mathcal{H}_0^*(\mathcal{C},\mathbb{T})$.

\begin{theorem} \label{thm2} Suppose that the function $f=h+\bar{g}$  belongs to $\mathcal{H}_0^*(\mathcal{C},\overline{\mathbb{D}})$ or $\mathcal{H}_0^*(\mathcal{C},\mathbb{T})$,
 where $h$ and $g$ are given by \eqref{equ1}.
Then $$|a_n|\leq \frac{n+1}{2}~\mbox{ and }~  |b_n|\leq\frac{n-1}{2}~\mbox{ for all $n \geq 2$.}
$$
The bounds are attained for the harmonic half-plane function $L$.
\end{theorem}

\bpf
The proof is similar to that of Theorem \ref{thm1}. We just give a sketch of proof.
Suppose that $f=h+\bar{g} \in \mathcal{H}_0^*(\mathcal{C},\overline{\mathbb{D}}),$ we have $f^\epsilon(z)=z+\sum_{n=2}^{\infty} \epsilon_{n} z^{n} \in \mathcal{C}$ for some $\epsilon\in\overline{\D}.$
It is a basic fact that $|\epsilon_{n}| \leq 1$ for all $n \geq 2$.
The desired result follows if the function $\frac{1+z}{(1-z)^3}$ in the proof of Theorem \ref{thm1} is replaced by
$\frac{1}{(1-z)^2}$ in two places.
In addition, we have $\mathcal{H}_0^*(\mathcal{C},\mathbb{T})\subset\mathcal{H}_0^*(\mathcal{C},\overline{\mathbb{D}})$ and $L\in\mathcal{H}_0^*(\mathcal{C},\mathbb{T})$.
\epf

\bcor \label{cor2}
Suppose that the function $f(z)=h(z)+\overline{g(z)}=z+\sum_{n=2}^{\infty} a_n z^n+\overline{\sum_{n=1}^{\infty} b_n z^n}$
belongs to $\mathcal{A}(\mathcal{H}_0^*(\mathcal{C},\overline{\mathbb{D}}))$ or $\mathcal{A}(\mathcal{H}_0^*(\mathcal{C},\mathbb{T}))$.
Then
$\left|a_{n}\right|<n$  and $\left|b_{n}\right|<n$ for all $n \geq 2$.
\ecor

For coefficient estimates on functions in classes $\mathcal{H}_0(\mathcal{S},\mathbb{T})$ and $\mathcal{H}_0(\mathcal{C},\mathbb{T})$,
please refer to  Theorem~8.1 and Propositions~8.2 in \cite{HM}, respectively.

\section{Growth, covering and distortion for subclasses of $\mathcal{H}_0^*(\mathcal{S},\overline{\mathbb{D}})$}  \label{sec4}

In this section, we will build the growth, covering and distortion theorems for functions in some subclasses of $\mathcal{A}(\mathcal{H}_0^*(\mathcal{S},\overline{\mathbb{D}}))$
along the similar proofs as  in \cite[Theorems~2,~3~and~4]{PS}. These results are based on the following lemma concerning affine and linear invariant family (ALIF).
We say that $\mathcal{F}\subseteq\mathcal{H}$ is ALIF if it is both AIF and LIF.
Here, a family $\mathcal{F}\subseteq\mathcal{H}$  is said to be an {\it affine invariant family} (AIF) if for each $f\in\mathcal{F}$,
\bes
A_c(f(z))=\frac{f(z)+c\overline{f(z)}}{1+c g'(0)}\in\mathcal{F}
\quad \forall ~c\in\mathbb{D},
\ees
and $\mathcal{F}\subseteq\mathcal{H}$ is called a {\it linear invariant family} (LIF) if for each $f\in\mathcal{F}$,
\bes
K_\varphi(f(z))=\frac{f(\varphi(z))-f(\varphi(0))}{\varphi'(0)h'(\varphi(0))}\in\mathcal{F}
\quad \forall ~\varphi\in{\rm Aut}(\mathbb{D}),
\ees
 where ${\rm Aut}(\mathbb{D})$ denotes the conformal automorphism of the unit disk $\mathbb{D}$.



\blem \label{lem1}
The family $\mathcal{A}(\mathcal{F})$ is an {\rm ALIF}, where $\mathcal{F}$ is chosen by anyone of six classes:   $\mathcal{H}_0(\mathcal{S},\mathbb{T})$,  $\mathcal{H}_0^*(\mathcal{S},\mathbb{T})$, $\mathcal{H}_0^*(\mathcal{S},\overline{\mathbb{D}})$, $\mathcal{H}_0(\mathcal{C},\mathbb{T})$, $\mathcal{H}_0^*(\mathcal{C},\mathbb{T})$ and $\mathcal{H}_0^*(\mathcal{C},\overline{\mathbb{D}})$.
\elem

\bpf
We will only prove that $\mathcal{A}(\mathcal{H}_0^*(\mathcal{S},\overline{\mathbb{D}}))$ is an ALIF, since the proofs for the other families are similar.
Let $f=h+\overline{g}\in\mathcal{A}(\mathcal{H}_0^*(\mathcal{S},\overline{\mathbb{D}}))$  be given.
Then there exists a constant $b_1\in\mathbb{D}$ and a function $f_0=h_0+\overline{g_0}\in\mathcal{H}_0^*(\mathcal{S},\overline{\mathbb{D}})$
such that $f=f_0+\overline{b_1f_0}$, where $h_0+\epsilon g_0\in\mathcal{S}$ for some $\epsilon\in\overline{\mathbb{D}}$.

The class $\mathcal{A}(\mathcal{H}_0^*(\mathcal{S},\overline{\mathbb{D}}))$ is an ALF, since
\bes
A_c(f(z))=f_0(z)+\frac{c+\overline{b_1}}{1+b_1 c}\overline{f_0(z)}\in\mathcal{A}(\mathcal{H}_0^*(\mathcal{S},\overline{\mathbb{D}}))
\quad \forall ~c\in\mathbb{D}.
\ees
To complete the proof, we only need to prove  that the class $\mathcal{A}(\mathcal{H}_0^*(\mathcal{S},\overline{\mathbb{D}}))$ is a LIF.
Let us consider the mapping $K_\varphi(f(z)):=H(z)+\overline{G(z)}$, where $\varphi\in{\rm Aut}(\mathbb{D})$.
Direct computation shows that $K_\varphi(f(z))=F_0(z)+\overline{B_1F_0(z)}$, where $F_0$ and $B_1$ are defined by
$$F_0(z)=H_0(z)+\overline{G_0(z)}=\frac{H(z)-\overline{B_1}G(z)}{1-|B_1|^2}
+\overline{\left(\frac{G(z)-B_1H(z)}{1-|B_1|^2}\right)}\quad \text{and}\quad
B_1=\frac{\varphi'(0)g'(\varphi(0))}{\overline{\varphi'(0)h'(\varphi(0))}},
$$
respectively. Clearly, $F_0\in\mathcal{H}_0$ and $B_1\in\mathbb{D}$.

Next we will prove $F_0\in\mathcal{H}_0^*(\mathcal{S},\overline{\mathbb{D}})$.
Let $$\rho=\frac{\overline{\varphi'(0)}}{\varphi'(0)}
\frac{\epsilon\overline{h'_0(\varphi(0))}+\overline{g'_0(\varphi(0))}}{h'_0(\varphi(0))+\epsilon g'_0(\varphi(0))}
=\frac{\overline{\varphi'(0)}}{\varphi'(0)}\frac{\overline{h'_0(\varphi'(0))}}{h'_0(\varphi'(0))}
\frac{\epsilon+\overline{\omega_{f_0}(\varphi(0))}}{1+\epsilon\omega_{f_0}(\varphi(0))},
$$
where $\omega_{f_0}(\varphi(0))=g'_0(\varphi(0))/h'_0(\varphi(0))$.
Clearly, $\rho\in\overline{\mathbb{D}}$.
Straightforward calculation yields that $H_0(z)+\rho G_0(z)=K_{\varphi}((h_0+\epsilon g_0)(z))$.
Since $h_0+\epsilon g_0\in\mathcal{S}$ and the family $\mathcal{S}$ is a LIF, we know $H_0+\rho G_0\in\mathcal{S}$,
which means that $F_0\in\mathcal{H}_0^*(\mathcal{S},\overline{\mathbb{D}})$.
Note that $\epsilon\in\mathbb{T}$ iff $\rho\in\mathbb{T}$. We complete the proof.
\epf

It is a basic fact that the family $\mathcal{S}_H=\mathcal{A}(\mathcal{S}_H^0)$ is an ALIF.
Motivated by the following growth and covering theorem for $\mathcal{S}_H^0$, we can obtain our corresponding theorems for  our subclasses of $\mathcal{H}_0^*(\mathcal{S},\overline{\mathbb{D}})$.
Furthermore, we need a distortion theorem concerning ALIF to build our distortion theorems for functions in restricted classes.
Now we will list the two results. For convenience, let $\alpha(\mathcal{F})=\ds \sup_{ f\in\mathcal{F}\subseteq \mathcal{H}} |a_2(f)|$.

\begin{Lem} {\rm{(\cite[Theorem,~p.\,97]{dur2004})}} \label{LemA}
Every function $f \in \mathcal{S}_H^0$ satisfies the inequalities
$$
\frac{1}{2 \alpha}\left[1-\left(\frac{1-r}{1+r}\right)^\alpha\right] \leq|f(z)| \leq \frac{1}{2 \alpha}\left[\left(\frac{1+r}{1-r}\right)^\alpha-1\right], \quad r=|z|<1 .
$$
In particular, the range of each function $f \in \mathcal{S}_H^0$ contains the disk $|w|<\frac{1}{2 \alpha}$. Here, $\alpha=\alpha(\mathcal{S}_H)$.
\end{Lem}

\begin{Lem} {\rm{(\cite[Theorem~1]{GE})}} \label{LemB}
Let $\mathcal{F}$ be an ALIF. Let $\mathcal{F}_0=\{f=h+\overline{g}\in\mathcal{F}:~g'(0)=0\}$ and $\alpha=\alpha(\mathcal{F}_0)$.
Suppose that $f \in \mathcal{F}$ with $b_1=f_{\bar{z}}(0)$.
Then the Jacobian $J_f$ of the mapping $f$ with any $z \in \mathbb{D}$ satisfies the bounds
$$
\left(1-\left|b_1\right|^2\right) \frac{(1-|z|)^{2 \alpha-2}}{(1+|z|)^{2 \alpha+2}} \leq J_f(z) \leq\left(1-\left|b_1\right|^2\right) \frac{(1+|z|)^{2 \alpha-2}}{(1-|z|)^{2 \alpha+2}}.
$$
Moreover, for any $z$ with $0<|z|=r<1$, the inequalities
$$
\left|h^{\prime}(z)\right| \leq\left(1+r\left|b_1\right|\right) \frac{(1+r)^{\alpha-3 / 2}}{(1-r)^{\alpha+3 / 2}} \quad \text { and } \quad\left|g^{\prime}(z)\right| \leq\left(r+\left|b_1\right|\right) \frac{(1+r)^{\alpha-3 / 2}}{(1-r)^{\alpha+3 / 2}}
$$
hold. These bounds are sharp for the class of univalent close-to-convex harmonic functions.
The equality is attained for the close-to-convex functions $f(z)=K+\overline{b_1 K}$, where $K$ is the harmonic Koebe function.
\end{Lem}

It is pointed out in \cite{dur2004} that Lemma~A 
remains valid for any  ALIF ~$\mathcal{M} \subseteq \mathcal{S}_H$ for which $\alpha =\alpha (\mathcal{M})$. 
Therefore, if the classes  $\mathcal{S}_H$ and $\mathcal{S}_H^0$  in Lemma~A 
are replaced by $\mathcal{S}_H\cap\mathcal{A}(\mathcal{H}_0^*(\mathcal{S},\overline{\mathbb{D}}))$ (resp. $\mathcal{S}_H\cap\mathcal{A}(\mathcal{H}_0^*(\mathcal{S},\mathbb{T}))$) and $\mathcal{S}_H^0\cap\mathcal{H}_0^*(\mathcal{S},\overline{\mathbb{D}})$ (resp. $\mathcal{S}_H^0\cap\mathcal{H}_0^*(\mathcal{S},\mathbb{T})$), respectively, then we can get corresponding results of Theorem \ref{thm3} by similar proofs.
The proofs are on  the basis of Lemma \ref{lem1} and a fact that $\mathcal{F}_1\cap\mathcal{F}_2~(\neq\emptyset)$ is an ALIF when both  $\mathcal{F}_1$ and $\mathcal{F}_2$ are ALIF.
Moreover, we use $\alpha(\mathcal{A}(\mathcal{H}_0^*(\mathcal{S},\overline{\mathbb{D}})))=\alpha(\mathcal{A}(\mathcal{H}_0^*(\mathcal{S},\mathbb{T})))=3$ by Corollary \ref{cor1} and $K\in\mathcal{H}_0^*(\mathcal{S},\mathbb{T})\cap\mathcal{S}_H^0$.
Similarly, we can get Theorem \ref{thm4} by applying Corollary \ref{cor2}, Lemmas \ref{lem1} and A.
We mention that $\mathcal{H}_0^*(\mathcal{C},\overline{\mathbb{D}})\subset\mathcal{S}_H^0$ and thus, $\mathcal{A}(\mathcal{H}_0^*(\mathcal{C},\overline{\mathbb{D}}))\subset\mathcal{A}(\mathcal{S}_H^0)=\mathcal{S}_H$.

\begin{theorem} {\rm (Growth and covering theorem I)}\label{thm3} Every function $f \in \mathcal{S}_{H}^{0}\cap\mathcal{H}_0^*(\mathcal{S},\overline{\mathbb{D}})$
(resp. $f\in\mathcal{S}_{H}^{0}\cap\mathcal{H}_0^*(\mathcal{S},\mathbb{T})$) satisfies the inequalities
$$
\frac{1}{6}\left[1-\left(\frac{1-r}{1+r}\right)^{3}\right] \leq|f(z)| \leq \frac{1}{6}\left[\left(\frac{1+r}{1-r}\right)^{3}-1\right], \quad r=|z|<1.
$$
In particular, the range of each function $f \in \mathcal{S}_{H}^{0}\cap\mathcal{H}_0^*(\mathcal{S},\overline{\mathbb{D}})$
(resp. $f\in\mathcal{S}_{H}^{0}\cap\mathcal{H}_0^*(\mathcal{S},\mathbb{T})$) contains the disk $|w|<\frac{1}{6}.$
The above inequalities are sharp and the equality is attained for the harmonic Koebe function $K$.
\end{theorem}

\begin{theorem} {\rm (Growth and covering theorem II)}\label{thm4} Every function $f \in \mathcal{H}_0^*(\mathcal{C},\overline{\mathbb{D}})$ (resp. $f\in\mathcal{H}_0^*(\mathcal{C},\mathbb{T})$) satisfies the inequalities
$$
\frac{1}{4}\left[1-\left(\frac{1-r}{1+r}\right)^{2}\right] \leq|f(z)| \leq \frac{1}{4}\left[\left(\frac{1+r}{1-r}\right)^{2}-1\right], \quad r=|z|<1.
$$
In particular, the range of each function $f \in \mathcal{H}_0^*(\mathcal{C},\overline{\mathbb{D}})$ (resp. $f\in\mathcal{H}_0^*(\mathcal{C},\mathbb{T})$) contains the disk $|w|<\frac{1}{4} .$
The above inequalities are sharp and the equality is attained for the harmonic half-plane function $L$.
\end{theorem}

For the growth and covering theorems for functions in $\mathcal{H}_0(\mathcal{S},\mathbb{T})$ and $\mathcal{H}_0(\mathcal{C},\mathbb{T})$, 
please refer to \cite[Proposition~8.3]{HM}.

Next by applying Theorems \ref{thm1} and \ref{thm2}, and Lemmas \ref{lem1} and B, 
we can get the following two results easily. Note that
$$\alpha(\mathcal{H}_0^*(\mathcal{S},\overline{\mathbb{D}}))=\alpha(\mathcal{H}_0^*(\mathcal{S},\mathbb{T}))=5/2
$$
and
$$\alpha(\mathcal{H}_0^*(\mathcal{C},\overline{\mathbb{D}}))=\alpha(\mathcal{H}_0^*(\mathcal{C},\mathbb{T}))=3/2.
$$

\begin{theorem} \label{thm5} {\rm (Distortion theorem I)}
Suppose that $f=h+\overline{g}$ belongs to $\mathcal{A}(\mathcal{H}_0^*(\mathcal{S},\overline{\mathbb{D}}))$  or $\mathcal{A}(\mathcal{H}_0^*(\mathcal{S},\mathbb{T}))$ with $b_{1}=f_{\bar{z}}(0).$
Then the Jacobian $J_{f}$ of the mapping $f$ with any $z \in \mathbb{D}$ satisfies the bounds
$$\left(1-\left|b_{1}\right|^{2}\right) \frac{(1-|z|)^{3}}{(1+|z|)^{7}} \leq J_{f}(z) \leq\left(1-\left|b_{1}\right|^{2}\right) \frac{(1+|z|)^{3}}{(1-|z|)^{7}}.
$$
Moreover, for any $z$ with $0<|z|=r<1,$ the inequalities
$$\left|h'(z)\right| \leq\left(1+r\left|b_{1}\right|\right) \frac{(1+r)}{(1-r)^{4}} \quad \text { and } \quad\left|g'(z)\right| \leq\left(r+\left|b_{1}\right|\right) \frac{(1+r)}{(1-r)^{4}}
$$
hold. These bounds are sharp. The equality is attained for the close-to-convex functions $f=K+\overline{b_{1} K}$, where $K$ is the harmonic Koebe function.
\end{theorem}

\begin{theorem} \label{thm6} {\rm (Distortion theorem II)}
Suppose that $f=h+\overline{g}$ belongs to $\mathcal{A}(\mathcal{H}_0^*(\mathcal{C},\overline{\mathbb{D}}))$ or $\mathcal{A}(\mathcal{H}_0^*(\mathcal{C},\mathbb{T}))$ with $b_{1}=f_{\bar{z}}(0) .$
Then the Jacobian $J_{f}$ of the mapping $f$ with any $z \in \mathbb{D}$ satisfies the bounds
$$\left(1-\left|b_{1}\right|^{2}\right) \frac{1-|z|}{(1+|z|)^{5}} \leq J_{f}(z) \leq\left(1-\left|b_{1}\right|^{2}\right) \frac{1+|z|}{(1-|z|)^5}.
$$
Moreover, for any $z$ with $0<|z|=r<1,$ the inequalities
$$\left|h'(z)\right| \leq \frac{1+r|b_{1}|}{(1-r)^3} \quad \text { and } \quad |g'(z)| \leq\frac{r+|b_1|}{(1-r)^3}
$$
hold. These bounds are sharp. The equality is attained for the convex functions $f=L+\overline{b_{1} L}$, where $L$ is the harmonic half-plane function.
\end{theorem}

Next we will build distortion theorems for functions in $\mathcal{A}(\mathcal{H}_0(\mathcal{S},\mathbb{T}))$ and $\mathcal{A}(\mathcal{H}_0(\mathcal{C},\mathbb{T}))$.
If we compare the last two theorems, one can easily see that there are some  differences in proofs.

\begin{theorem} \label{thm7} {\rm (Distortion theorem III)}
Let $f=h+\overline{g}\in\mathcal{A}(\mathcal{H}_0(\mathcal{S},\mathbb{T}))$ with $b_{1}=f_{\bar{z}}(0).$
Then the Jacobian $J_{f}$ of the mapping $f$ with any $z \in \mathbb{D}$ satisfies the bounds
$$\left(1-\left|b_{1}\right|^{2}\right) \frac{(1-|z|)^2}{(1+|z|)^6} \leq J_{f}(z) \leq\left(1-\left|b_{1}\right|^{2}\right) \frac{(1+|z|)^2}{(1-|z|)^6}.
$$
Moreover, for any $z$ with $0<|z|=r<1,$ the inequalities
$$\left|h'(z)\right| \leq\frac{1+r}{(1-r)^3} \quad \text { and } \quad\left|g'(z)\right| \leq \min\{1,r+|b_1|\}\frac{1+r}{(1-r)^3}
$$
hold. The former three  inequalities are sharp and the equality is attained for  functions $f=k+\overline{b_1k}$, where $k$ is the  Koebe function.
\end{theorem}

\bpf
It follows from \cite[Theorem~8.1]{HM} that $\alpha(\mathcal{H}_0(\mathcal{S},\mathbb{T}))=2$.
By using Lemmas \ref{lem1} and B, 
we can get the former two inequalities.

Next we consider the last two inequalities. Let $f=h+\overline{g}\in\mathcal{A}(\mathcal{H}_0(\mathcal{S},\mathbb{T}))$ with $b_{1}=f_{\bar{z}}(0).$
Then there exists a function $f_0=h_0+\overline{g_0}\in\mathcal{H}_0(\mathcal{S},\mathbb{T})$ such that $f=f_0+\overline{b_1f_0}$.
Clearly, we have $h=h_0+\overline{b_1}g_0$ and $g=b_1h_0+g_0$ in $\ID$.
It follows from \eqref{equ4} that $h_0+\epsilon g_0\in\mathcal{S}$ for every $\epsilon\in\overline{\mathbb{D}}$,
and thus $h\in\mathcal{S}$ and $h_0\in\mathcal{S}$.
Therefore, the third inequality follows from the classical distortion theorem for functions in $\mathcal{S}$.
For the last inequality, on one hand, we have $|g'|<|h'|$ in $\ID$ since $f$ is sense-preserving.
On the other hand, by Schwarz-Pick lemma, we can get
$$|g'(z)|=|b_1h'_0(z)+g'_0(z)|=|h'_0(z)| |b_1+\omega_{f_0}(z)|\leq(|b_1|+|z|)| h'_0(z)|,
$$
and thus, the desired result follows from the third inequality. We finish the proof.
\epf

Similarly, we can get the following result associated \eqref{equ4}, \cite[Proposition~8.2]{HM} with Lemmas \ref{lem1} and B. 

\begin{theorem} \label{thm8} {\rm (Distortion theorem IV)}
Let $f=h+\overline{g}\in\mathcal{A}(\mathcal{H}_0(\mathcal{C},\mathbb{T}))$ with $b_{1}=f_{\bar{z}}(0) .$
Then the Jacobian $J_{f}$ of the mapping $f$ with any $z \in \mathbb{D}$ satisfies the bounds
$$\left(1-\left|b_{1}\right|^{2}\right) \frac{1}{(1+|z|)^2} \leq J_{f}(z) \leq\left(1-\left|b_{1}\right|^{2}\right) \frac{1}{(1-|z|)^2}.
$$
Moreover, for any $z$ with $0<|z|=r<1,$ the inequalities
$$\left|h'(z)\right| \leq\frac{1}{(1-r)^2} \quad \text { and } \quad\left|g'(z)\right| \leq \min\{1,r+|b_1|\}\frac{1}{(1-r)^2}
$$
hold. The former three inequalities are sharp and the equality is attained for  functions $f=l+\overline{b_1l}$, where $l(z)=z/(1-z)$.
\end{theorem}

\section{Bohr inequalities for functions in  the subclasses of $\mathcal{H}_0^*(\mathcal{S},\overline{\mathbb{D}})$}  \label{sec5}
In this section, we will build some Bohr inequalities for functions subordinate to  functions in the subclasses of $\mathcal{H}_0^*(\mathcal{S},\overline{\mathbb{D}})$.
The obtained results are generalizations of the  earlier known result (\cite[Theorem~1]{abu}), which establishes  Bohr's inequality for the class of analytic functions subordinate to univalent functions.
Such result generates the classical Bohr's inequality  which asserts that
$\sum_{n=0}^{\infty} |a_n|r^n\leq 1$ holds for $|z|=r\leq1/3$ if $f(z)=\sum_{n=0}^{\infty} a_nz^n$ with $|f(z)|\leq1$ in $\ID$.
In fact, Bohr obtained it for $|z|\leq1/6$ (see \cite{boh}), which was sharpened to $|z|\leq1/3$ by Wiener, Riesz and Schur independently.
The constant $1/3$ is optimal and is called the Bohr radius.


It is necessary to recall the notion of subordination for harmonic mappings.
For two harmonic mappings $f=h+\overline{g}$ and $F=H+\overline{G}$ defined in $\ID$, $f$ is said to be subordinate to $F$,
written by $f\prec F$, if there is a  function $\omega$ analytic in $\ID$ with $\omega(0)=0$ and $|\omega(z)|<1$ such that $f(z)=F(\omega(z))$,
i.e., $h(z)=H(\omega(z))$ and $g(z)=G(\omega(z))$ (see \cite{sch}).
Clearly, $h\prec H$ and $g\prec G$ in the sense of subordination for analytic functions.
In particular, when $H$ is univalent in $\D$, then the relationship $h \prec H$ is equivalent
to the conditions $h(0)=H(0)$ and $h(\D)\subset H(\D)$.
Such a geometric property cannot be extended to subordination for harmonic mappings in general.
Throughout, we denote the class of all functions $f$ subordinate to a fixed function $F$
and  the distance from the point $F(0)$ to the boundary $\partial F(\D)$ of the domain $F(\D)$ by $S(F)$ and $\dist(F(0),\partial F(\D))$, respectively.

Below we will use the following result several times.

\begin{Lem}\label{LemC}
{\rm (\cite[Corollary~2.2]{AlkKayPon} and \cite[Lemma~1]{BD})}
Let $f(z)$ and $g(z)$ be two analytic functions in $\ID$ with the Taylor series expansions
$f(z)=\sum_{k = 0}^\infty a_kz^k$ and $g(z)=\sum_{k = 0}^\infty b_kz^k$ for $z\in\ID$. If $f(z)\prec g(z)$,  then
$$\sum_{k=0}^\infty |a_k| r^k \le \sum_{k=0}^\infty |b_k| r^k \  \ \text{for all}\  \ r \le \frac{1}{3}.
$$
The number $1/3$ is optimal.
\end{Lem}

The fact that $1/3$ is sharp has been shown in \cite{PVW} and the last inequality holds for a general class of quasi-subordination \cite{AlkKayPon}.

\begin{theorem} \label{thm9}
 Let $f(z)=h(z)+\overline{g(z)}=\sum_{n=1}^{\infty} a_n z^n+\overline{\sum_{n=1}^{\infty} b_n z^n}\in S(F)$ and $F\in\mathcal{S}_H^0\cap\mathcal{H}_0^*(\mathcal{S},\overline{\mathbb{D}})$
 or $F\in\mathcal{S}_H^0\cap\mathcal{H}_0^*(\mathcal{S},\mathbb{T})$.
Then
$$
\sum_{n=1}^{\infty}\left|a_n\right| r^n+\sum_{n=1}^{\infty}\left|b_n\right| r^n \leq \dist(F(0),\partial F(\D))
$$
for $|z|=r \leq r_1=0.115013\ldots$, where $r_1$ is the unique root in $(0,1)$ of $3r^{3}-3 r^2+9r-1=0$.
The number $r_1$ is optimal.
\end{theorem}

\bpf
We only consider the case $F\in\mathcal{S}_H^0\cap\mathcal{H}_0^*(\mathcal{S},\overline{\mathbb{D}})$,
 and thus, we can assume that $F=H+\overline{G}$,
where $H(z)=z+\sum_{n=2}^{\infty} A_n z^n$ and $G(z)=\sum_{n=2}^{\infty} B_n z^n$.
The notation of $f\in S(F)$ implies that $h\prec H$ and $g\prec G$.
By Lemma~C, 
we have
$$\sum_{n=1}^{\infty} |a_n| r^n\leq r+\sum_{n=2}^{\infty} |A_n| r^n\quad \text{and}
\quad \sum_{n=1}^{\infty} |b_n| r^n\leq \sum_{n=2}^{\infty} |B_n| r^n
$$
for $|z|=r\leq1/3$.
Since $F\in\mathcal{S}_H^0\cap\mathcal{H}_0^*(\mathcal{S},\overline{\mathbb{D}})$, in view of Theorem \ref{thm3}, we have
$$\dist(F(0),\partial F(\D))\geq1/6.
$$
To complete the proof, it is enough to show the sharp inequality
$$r+\sum_{n=2}^{\infty} |A_n| r^n+\sum_{n=2}^{\infty} |B_n| r^n \leq \frac{1}{6}
$$
holds for $|z|=r\leq r_1$.

It follows from Theorem \ref{thm1} that
\begin{align*}
 r+\sum_{n=2}^{\infty} |A_n| r^n+\sum_{n=2}^{\infty} |B_n| r^n
 \leq & \sum_{n=1}^{\infty} \frac{(n+1)(2 n+1)}{6} r^n+\sum_{n=2}^{\infty} \frac{(n-1)(2 n-1)}{6} r^n \\
 = &\frac{r-\frac{1}{2}r^2+\frac{1}{6}r^3}{(1-r)^3}+\frac{\frac{1}{2}r^2+\frac{1}{6}r^3}{(1-r)^3} \\
 = &\frac{r+\frac{1}{3}r^3}{(1-r)^3}
\end{align*}
which is less than or equal to $1/6$ whenever $3r^{3}-3 r^2+9r-1\leq 0$. This gives the desired inequality for $|z|=r\leq r_1$,
and the harmonic Koebe function shows that the number $r_1$ is best possible.
Note that $K$ maps $\ID$ univalently onto the slit plane $\mathbb{C}\setminus\{w=u+iv:~u\leq-1/6,~v=0\}$ (see \cite{CS}).
We complete the proof.
\epf

Similarly, we can get the following result by  Lemma~C, 
Theorems \ref{thm2} and \ref{thm4}.
The sharpness can be seen from the harmonic half-plane mapping $L$.
We omit the proof.

\begin{theorem} \label{thm10}
 Let $f(z)=h(z)+\overline{g(z)}=\sum_{n=1}^{\infty} a_n z^n+\overline{\sum_{n=1}^{\infty} b_n z^n}\in S(F)$ and $F\in\mathcal{H}_0^*(\mathcal{C},\overline{\mathbb{D}})$
 or $F\in\mathcal{H}_0^*(\mathcal{C},\mathbb{T})$.
Then
$$
\sum_{n=1}^{\infty}\left|a_n\right| r^n+\sum_{n=1}^{\infty}\left|b_n\right| r^n \leq \dist(F(0),\partial F(\D))
$$
for $|z|=r \leq r_{2}=3-2\sqrt{2}\approx0.171573$.
The number $r_2$  is optimal.
\end{theorem}

The following result can be easily got by applying Lemma~C, 
and Theorem~8.1,~Propositions~8.2 and 8.3 in \cite{HM}.

\begin{theorem} \label{thm11}
Let $f(z)=h(z)+\overline{g(z)}=\sum_{n=1}^{\infty} a_n z^n+\overline{\sum_{n=1}^{\infty} b_n z^n}\in S(F)$ and
$F\in\mathcal{H}_0(\mathcal{S},\mathbb{T})$  (resp. $F\in\mathcal{H}_0(\mathcal{C},\mathbb{T})$).
Then
$$
\sum_{n=1}^{\infty}\left|a_n\right| r^n+\sum_{n=1}^{\infty}\left|b_n\right| r^n \leq \dist(F(0),\partial F(\D))
$$
for $|z|=r\leq 3-2\sqrt{2}$  $($resp. $|z|=r\leq1/3$ $)$.
The number $3-2\sqrt{2}$ $($resp. $1/3$ $)$ is optimal.
\end{theorem}

\subsection*{Conflict of Interests}
The authors declare that there is no conflict of interests regarding the publication of this paper.

\subsection*{Acknowledgments}
The research was partly supported by 
the Hunan Provincial Natural Science Foundation of China (No. 2021JJ30057),
the Science and Technology Plan Project of Hunan Province (No. 2016TP1020) and the Application-Oriented Characterized Disciplines, Double First-Class University Project of Hunan Province (Xiangjiaotong [2018]469).
The authors would like to thank   the referees for their valuable comments for improving this paper.


\begin{thebibliography}{30}

\bibitem{abu} Y. Abu Muhanna,
Bohr's phenomenon in subordination and bounded harmonic classes,
\emph{Complex Variables Theory Appl.} \textbf{55}(11) (2010), 1071--1078.


%



\bibitem{AlkKayPon} S. A. Alkhaleefah, I. R. Kayumov and S. Ponnusamy,
On the Bohr inequality with a fixed zero coefficient,
\emph{Proc. Amer. Math. Soc.} \textbf{147}(12) (2019), 5263--5274.

\bibitem{ABLS} M. Aydogan, D. Bshouty, A. Lyzzaik and F. M. Sakar, On the shears of univalent harmonic mappings,
\emph{Complex Anal. Oper. Theory} \textbf{13} (2019), 2853--2862.


\bibitem{BD} B. Bhowmik and N. Das, Bohr phenomenon for subordinating families of certain univalent functions,
\emph{J. Math. Anal. Appl.} \textbf{462}(2) (2018), 1087--1098.

%
%

\bibitem{boh} H. Bohr,
A theorem concerning power series,
\emph{Proc. London Math. Soc.} \textbf{13}(2) (1914), 1--5.

\bibitem{CH} M. Chuaqui and R. Hern\'{a}ndez, Univalent harmonic mappings and linearly connected domains,
\emph{J. Math. Anal. Appl.} {\bf 332}(2) (2007), 1189--1194.

\bibitem{CS} J. G. Clunie and T. Sheil-Small,
Harmonic univalent functions,
\emph{Ann. Acad. Sci. Fenn. Ser. A. I.} {\bf 9} (1984), 3--25.

\bibitem{deB} L. de Branges, A proof of the Bieberbach conjecture, \emph{Acta Math.} {\bf 154}(1-2) (1985), 137--152.


%


\bibitem{dur2004} P. Duren, Harmonic Mappings in the Plane,
\emph{Cambridge University Press}, New York, 2004.




 \bibitem{GE} S. Yu. Graf and O. R. Eyelangoli, Differential inequalities in linear- and affine- invariant families of harmonic mappings,
\emph{Russian Math.. (Iz. VUZ)} {\bf 54}(10) (2016), 60--62.

%


\bibitem{HM} R. Hern\'{a}ndez and  M. J. Mart\'{\i}n, Stable geometric properties of analytic and harmonic mappings,
\emph{Math. Proc. Camb. Phil. Soc.} {\bf 155} (2013), 343--359.


%

\bibitem{KPX} I. R. Kayumov, S. Ponnusamy and  L. A. Xuan,  Rotations of convex harmonic univalent mappings,
\emph{Bull. Sci. Math} \textbf{155} (2019), 1--9.

\bibitem{lew} H. Lewy,
On the non-vanishing of the Jacobian in certain one-to-one mappings,
\emph{Bull. Amer. Math. Soc.} {\bf 42} (1936), 689--692.


\bibitem{LP2019} G. Liu and S. Ponnusamy, Harmonic pre-Schwarzian and its application,
\emph{Bull. Sci. Math.} \textbf{152}(3) (2019),  150--168.

%



\bibitem{pom} C. Pommerenke, \emph{Univalent Functions}, Vanden-hoek and Ruprecht, Gottingen, 1975.
%
%
%

\bibitem{PS} S. Ponnusamy and A. Sairam Kaliraj,
On the coefficient conjecture of Clunie and Sheil-Small on univalent harmonic mappings,
\emph{Proc. Indian Acad. Sci.} {\bf 125}(3) (2015), 277--290.

\bibitem{PVW} S. Ponnusamy, R. Vijayakumar and K.-J. Wirths,
{Improved Bohr's phenomenon in quasi-subordination classes}
\emph{J. Math. Anal. and Appl.} \textbf{506}(1) (2022), 10 pages, Article 125645.

\bibitem{sch} L. E. Schaubroeck, Subordination of planar harmonic functions,
\emph{Complex Variables} {\bf 41} (2000), 163--178.

\bibitem{she} T. Sheil-Small, Constants for planar harmonic mappings,
\emph{J. London Math. Soc.} \textbf{42} (1990), 237--248.

\bibitem{suf} T. J. Suffridge, Some remarks on convex maps of the unit disk,
\emph{Duke Math.} {\bf 37} (1970), 775--777.

\bibitem{WLZ} X. Wang, X. Liang and Y. Zhang,
Precise coefficient estimates for close-to-convex harmonic univalent mappings,
\emph{J. Math. Anal. Appl.} \textbf{263} (2001), 501--509.
\end{thebibliography}
\end{document}